\documentclass[12pt]{amsart}
\usepackage[english]{babel}
\usepackage{amsmath}
\usepackage{amsthm}
\usepackage{amssymb}
\usepackage{mathrsfs}
\usepackage{enumerate}
\usepackage[notcite, final, notref]{showkeys}
\usepackage{dsfont}
\usepackage[dvips]{color}
\usepackage[colorlinks,linkcolor=red,anchorcolor=green,citecolor=blue]{hyperref}
\usepackage{xcolor}
\newtheorem{theorem}{Theorem}[section]

\newtheorem{lemma}[theorem]{Lemma}

\theoremstyle{definition}
\newtheorem{remark}[theorem]{Remark}

\def\mmA{\mathbb{A}^2_{\alpha}}
\def\D{\bf D}

\def\H{\mathbb{H}^2}

\def\C{\bf{C}}

\def\D{\bf D}

\title[]{\bf $n$-Best Kernel Approximation in Reproducing Kernel Hilbert Spaces}

\author[T. Qian]{Tao Qian}
\thanks{Email: tqian@must.edu.mo~~~ORCID: 0000-0002-8780-9958\\
Macau University of Science and Technology, (Macau Center for Mathematical Sciences), Macau, China\\
Funded by The Science and Technology Development Fund, Macau SAR (File no. 0123/2018/A3)}

\begin{document}
\maketitle
\begin{abstract}
By making a seminal use of the maximum modulus principle of holomorphic functions we prove existence of $n$-best kernel approximation for a wide class of reproducing kernel Hilbert spaces of holomorphic functions in the unit disc, and for the corresponding class of Bochner type spaces of stochastic processes. This study thus generalizes the classical result of $n$-best rational approximation for the Hardy space and a recent result of $n$-best kernel approximation for the weighted Bergman spaces of the unit disc. The type of approximations have significant applications to signal and image processing and system identification, as well as to numerical solutions of the classical and the stochastic type integral and differential equations.
\end{abstract}

\bigskip

\date{today}

\bigskip


\section{Introduction}

A main form of application of mathematical analysis is approximation by basic functions of the underlying space. Various forms and topics of polynomial and rational approximations have been studied, including convergence models, capacity and rates, existence and uniqueness of best approximation, as well as algorithms, etc.  See for instance the selected list of the literature \cite{J.L.W1,1963che,1964CL,1974S,1967C,1970L,1988SK,1989KD,1994BY,LB1,LB2,2001D,2003zhou,
2004BCL,2018,Seip} and the references therein. The present study will concentrate in approximation of reproducing kernel Hilbert spaces (RKHSs) of complex holomorphic functions, the latter being related to $Z$-transforms of system transfer functions. In RKHSs the most natural basic functions are the parameterized reproducing kernels. Tasks of signal and image processing are based on effective reconstruction of a given signal or image. To measure reconstruction efficiency, among the most commonly used, there are two dual models. One is, for a previously given $\epsilon$, to determine the smallest integer $n$ such that the difference, measured in the underlying space norm, between the given function and some $n$-linear combination of the basic functions is already dominated by $\epsilon.$  The second model is, for a given resource limitation represented by a natural integer $n,$ to find an $n$-tuple of parameterized basic functions and an $n$-tuple of coefficients such that the $n$-linear combination that they compose gives rise to the best possible approximation to the given function. The second model is abbreviated as $n$-\emph{best approximation} or more briefly $n$-\emph{best problem}. The present paper restricts to study the second model, and only the existence part of the $n$-best solutions. The corresponding algorithm part for actually finding one or all $n$-best solutions, as a consequence of the technical results of the existence proof, will be separately studied.\\

We will be based on the general concept reproducing kernel Hilbert space (RKHS). In 1930, to study the partial differential equation
\[\frac{\partial^2u}{\partial x^2}+\frac{\partial^2u}{\partial y^2}+\alpha (x,y)\frac{\partial u}{\partial x}+ \beta (x,y)\frac{\partial u}{\partial y}+
\gamma (x,y)u=0,\]
where $\alpha (x,y), \beta (x,y), \gamma (x,y)\in C^2(\Omega), \Omega$ is a bounded region, $\alpha (x,y), \beta (x,y), \gamma (x,y)\in C^2(\Omega), \Omega$ are all real-analytic functions, S. Bergman proposed the reproducing kernel concept and gave the related formulas. The theoretical study of reproducing kernel may be divided into two stages. Of which the first started from J. Mercer (\cite{Mer}), at the beginning of the 20th century, who in his studies of integral equations brought up the concept positive definite kernel:
 \[ \sum_{i,j=1}^n K(y_i,y_j)\xi_i\xi_j\ge 0.\]
The second stage is the development by E. H. Moore (\cite{Moo}), around 1930's, who proved that every positive Hermitian matrix induces a Hilbert space that has a kernel function $K(x,y)$ enjoying the property
\[ f(y)=\langle f, K(\cdot, y)\rangle.\]
The same phenomenon was observed by S. Bochner in the convolution kernel form connected to Fourier theory (\cite{B1}). Around 1940's the most popular reproducing kernels were the Bergman type ones. Bergman developed the idea of S. Zaremba (\cite{Zaremba}) to solve boundary value problems by using reproducing kernels, showing that reproducing kernels are effective tools to solve elliptic boundary value problems.\\

Combined with various science and engineering objects there have developed new theories and algorithms, including signal processing (\cite{AD}), system identification (\cite{AN,MQ,MQQ,QWxy}), scholastics processing (\cite{HH,QianSAFD,QQD,YQCZQ}), estimation theory (\cite{Gu,PKJ}), wavelet transform (\cite{SS}), reproducing kernel particle method (\cite{HLC,GLDC,GL,HM,CYWL,VC}), the moving least-square reproducing kernel method (\cite{LL,LLT}), multi-scale reproducing kernel particle methods (\cite{LCUC}), etc., with ample applications. \\

F. M. Larkin (\cite{Larkin}) and M. M. Chawla (\cite{Chawla}) studied the approximation aspect. Formats of interpolation with reproducing kernel have been used to numerical solutions of partial differential equations and integral equations.  The latest reproducing kernel approximation methodology, called adaptive Fourier decomposition (AFD), uses the maximal energy extraction principle similar to what is in greedy algorithm or matching pursuit (\cite{MZ,Te}). Greedy algorithms are based on general Hilbert space theory with a dictionary (\cite{Te2}).
The AFD methods, originated from analytic positive frequency decomposition, validate  attainability of the best suited parameters. Technically, AFDs are based on delicate mathematical analysis, and, in particular, allow repeating selection of parameters through defining multiple kernels, when necessary. The technical treatment is a blend of functional analysis, complex analysis and harmonic analytic. Recent studies on approximation in Hardy spaces, including the latest $n$-best and stochastic AFD approximations, may be found in a sequence of articles \cite{QWa,Q2D,ACQS1,ACQS2,CS,CP,Q222,WaQ,QianSAFD,QQD}. Some early studies of adaptive Fourier decomposition in Bergman and weighted Bergman spaces are given in \cite{qu2018,qu2019}. Celebrating results for the $n$-best type approximation in weighted Bergman spaces are presented in \cite{QQLZ}.
\\

Recent studies given by Ball et al.'s papers (\cite{Ball1,Ball2}) show that approximations of Hilbert spaces of holomorphic functions have intimate connections with system identification and in particular with time-variant linear systems. In the Hardy space case there holds the Sz.-Nagy-Foias model theory for $C_{\cdot 0}$ contraction operators. The model theory combined with the Burling-Lax theorem addresses a correspondence between any two of the four kinds of objects: shift invariant subspaces, operator-valued inner functions, conservative discrete-time input/state/output linear system, and $C_{\cdot 0}$ Hilbert-space contraction operators. The studies of \cite{Ball1} and \cite{Ball2} extend such correspondence to weighted Bergman and weighted Hardy spaces. Under such frame work and via the $Z$-transform of the system, $n$-best approximation in each of the mentioned spaces determines the optimal shift invariant subspaces for effective and efficient system identification.\\

 To complete the introduction for science and engineering motivations and involvements, we at last, but not least, mention that there have been forerunner but recent developments on the stochastic $n$-best model as cited in \cite{QianSAFD,QQD,YQCZQ}. Stochastic $n$-best approximation of our general setting will be given in \S 5. Stochastic AFD offers a new approach to stochastic processes, including solutions of stochastic partial differential equations. It, in particular, stands as an alternative method to the Karhunen-Lo\'eve decomposition together with several advantages. Further developments along this direction are to be reported in separate and forthcoming papers.  \\

Next, we introduce the related preliminary knowledge and raise the $n$-best
problem in our general setting.
 The Hardy space is defined
\[ \mathbb{H}^{2}({\D})=\{f:{\D}\to {\C} \ :\ f \ {\rm is\ analytic\ in } \ {\D}\ {\rm and}\ \|f\|_{\H}=\sup_{0\leq r<1}\int_0^{2\pi} |f(r{\rm e}^{it})|^2<\infty\}.\]
As a fundamental result, functions in the Hardy space have non-tangential boundary limits a.e. on the unit circle $\partial{\D}$. The Hardy space is isometric to the function space ${\H}_{\partial{\D}}$ consisting of the non-tangential limiting functions. One of the alternative definitions of the space ${\H}_{\partial{\D}}$ is
\[ {\H}_{\partial{\D}}=\{f: \partial{\D}\to {\C}\ :\ f\in L^2(\partial {\D}), \quad
f({\rm e}^{it})=\sum_{k=0}^\infty c_k{\rm e}^{ikt},\quad \sum_{k=0}^\infty |c_k|^2<\infty\}.\]
${\H}_{\partial{\D}}$ is a closed subspace of the Hilbert space $L^2(\partial {\D})$ equipped with the inner product
\begin{eqnarray}\label{inner} \langle f,g\rangle=\frac{1}{2\pi}\int_0^{2\pi}f({\rm e}^{it})\overline{g({\rm e}^{it})}dt.\end{eqnarray}

For a given positive integer $n,$ an ordered pair of polynomials $(p,q)$ is called an $n$-\emph{admissible pair} if $p$ and $q$ are co-prime, $q\ne 0$ in ${\D},$ and both the degrees of $p$ and $q$ do not exceed $n.$ The famous $n$-\emph{best rational approximation problem} in $\mathbb{H}^{2}({\D})$ is as follows: For $f\in \mathbb{H}^{2}({\D}),$ find an $n$-admissible ordered pair $(\tilde{p},\tilde{q})$ such that
\begin{eqnarray}\label{problem}
\|f-\frac{\tilde{p}}{\tilde{q}}\|_{\mathbb{H}^2({\D})}=\inf \{\|f-\frac{p}{q}\|_{\mathbb{H}^2({\D})} \ :\ (p,q) \ {\rm is\ an} \ n{\mbox -admissible\ pair}\}.\end{eqnarray}

 Closely related to rational approximation there exist studies on what is called Takenaka-Malmquist (TM)
system, or rational orthogonal system:
\begin{eqnarray}\label{B}
\{B_{a_1a_2...a_n}(z)\}_{n=1}^\infty=\left\{\frac{\sqrt{1-|a_{n}|^{2}}}{1-\bar{a}_{n}z}
\prod_{k=1}^{n-1}\frac{z-a_{k}}{1-\bar{a}_{k}z}\right\}_{n=1}^\infty,\quad a_1,\cdots,a_n,\cdots, \in {\D},
\end{eqnarray}
where $B_{a_1a_2...a_n}(z)=e_{a_n}(z)\phi_{a_1a_2...a_{n-1}}(z),$ where $e_{a_n}$ is the normalized Szeg\"o kernel at $a_n,$
\[e_{a_n}(z)=\frac{\sqrt{1-|a_n|^2}}{1-\overline{a}_nz
},\]
and the \emph{canonical Blaschke product} with zeros $a_1,\cdots,a_{n-1},$ \begin{eqnarray}\label{phi}\phi_{a_1a_2...a_{n-1}}(z)=
\prod_{k=1}^{n-1}\frac{z-a_{k}}{1-\bar{a}_{k}z}.\end{eqnarray}
There further holds the relation, for $a\in \D,$
\begin{eqnarray}\label{e} e_a(z)=\frac{\sqrt{1-|a|^2}}{1-\overline{a}z}=\frac{k_a}{\|k_a\|},\end{eqnarray}
where
\[ k_a(z)=\frac{{1}}{1-\overline{a}z} \quad {\rm and}\quad \|k_a\|=\frac{1}{\sqrt{1-|a|^2}}\]
 are, respectively, the reproducing kernel of $\mathbb{H}^{2}({\D})$ and the normalizing constant making $\|e_a\|=1.$\\

 The above formulated $n$-best rational approximation problem (\ref{problem}) is, in essence,  equivalent to the following $n$-\emph{best Blaschke form approximation} problem (\cite{QWe,Q222}): Let $n$ be a given positive integer. If $f$ itself is not an $m$-Blaschke form for some $m<n,$ find a set of $n$ parameters $a_1,\cdots,a_n,$ all in ${\D},$  such that
\begin{eqnarray}\label{nbest} \|f-\sum_{k=1}^n\langle f,B_{a_1\cdots a_k}\rangle B_{a_1\cdots a_k}\|&=&\inf \{\|f-\sum_{k=1}^n\langle f,B_{b_1\cdots b_k}\rangle  B_{b_1\cdots b_k}\|\ :\ b_1,\cdots,b_n\ \in {\D}\}.\nonumber
\end{eqnarray}
 The formulation of the problem allows multiplicity of the zeros $a_k$'s.
A solution of (\ref{nbest}) will be referred as an $\lq\lq n$-best Blaschke form approximation"(\cite{Q222,QWe}). Regardless unimodular constants (see Lemma \ref{lemma2} below), an $n$-TM system is the Gram-Schmidt orthogonalization of a set of $n$  Szeg\"o kernels, or multiple kernels (formulated in (\ref{Hardy multiple kernel}) below) when the parameters are with multiplicities. See explanations in \S 2). The $n$-best Blaschke form approximation is, again, equivalent with the $n$-\emph{best kernel approximation} formulated as: Find $(a_1,\cdots,a_n)\in {\D}^n$ and $(c_1,\cdots,c_n)\in {\C}^n$ such that
\begin{eqnarray}\label{nbest kernel}\|f-\sum_{j=1}^n c_j\tilde{k}_{a_j} \|
=\inf \{\|f-\sum_{j=1}^n c'_j{k}_{a'_j}\|: a'_1,\cdots,a'_n\ \in {\D}\\ \nonumber
 {\rm are\ distinct}, {\rm and} c'_1,\cdots,c'_n\ \in \C\},
\end{eqnarray}
where $\tilde{k}_{a_j}$ are multiple kernels (see (\ref{Hardy multiple kernel})).\\

Considerable amount of studies have been devoted to the above problem with the three equivalent forms. See \cite{J.L.W1,J.L.W2,LB1,LB2,Q222,QWe,MQ,QQLZ}. Amongst, researchers have obtained several new proofs for existence of a solution to (\ref{problem}). The motivation of exploring new proofs of the existence, including that of the author himself's, would be at least two-folder: (i) The known existence proofs for the Hardy space case do not seem to be adaptable to prove existence of an $n$-best approximation in any non-Hardy spaces, including weighted Bergman spaces and weighted Hardy spaces. In fact, before the work \cite{QQLZ} whether there is a solution to the problem in any non-Hardy space was unknown; and (ii) On top of the existence, an ultimate algorithm of finding even one $n$-best solution of (\ref{problem}) has yet to be sought: The commonly adopted empirical algorithms are all local that cannot theoretically avoid the possibility of sinking into a local minimum. See \cite{LB2,LB1,QWa,Q2D,qu2018,QQ2} and the references therein. As an extension of the traditional Fourier method, both the n-best and the repeated one-by-one types have been found to have effective applications in signal and image processing and system identification (\cite{QWxy,QLyt,MQ}).\\

The question for $n$-best kernel approximation can be raised in general RKHSs, or even in Hilbert spaces with a dictionary. The recently published new proof of existence of $n$-best approximation in the Hardy space case (\cite{WQ2020}) uses the maximum modulus principle of holomorphic functions as a crucial technical trick. Inspired by this complex analysis method, through proving some necessary new pointwise estimations of the kernels of the involved zero- and Blaschke-weighted spaces, the study \cite{QQLZ} was able to prove existence of the $n$-best kernel approximation of all the weighted Bergman spaces $\mmA({\D}), -1<\alpha <\infty.$ \\

The study in the present paper is directly motivated by the recent new proofs of the existence on the Hardy space (\cite{WaQ}) and one on the weighted Bergman spaces (\cite{QQLZ}). We achieve a clever and concise proof for the existence of the $n$-best kernel approximation for a large class of reproducing kernel Hilbert spaces, that is in particular strictly larger than that of the Bergman ones, enclosing all the weighted Hardy spaces.\\

Based on further analysis of the orthogonalization projection operator $Q_{a_1a_2...a_k}$ and factorization of higher order generalized backwardd
 shift operators $Q_{a_1a_2...a_k}/\phi_{a_1\cdots{a_{m-1}}}$ (see \S 3), the present paper is able to avoid use of the pointwise estimations of the reproducing kernels of the involved zero spaces and the Blaschke weighted spaces. For RKHSs more general than the Bergman ones such kernel estimations may be impossible.  As a result, we are able to assert existence of the $n$-best problem for a class of RKHSs more general  than the weighted Hardy spaces. Precisely, we can declare existence of solutions of the $n$-best kernel approximation for all RKHSs of holomorphic functions in $\D$ that satisfy the following three conditions (see Theorem \ref{Theorem 1}):\\

\noindent (i) The reproducing kernel $K(z,w)$ enjoys \emph{the analyticity condition}: When $w\in \D$ is fixed, $K(z,w)$ is analytic for $z$ in a neighbourhood of the closed unit disc $\overline{\D},$ and, when $z\in \D$ is fixed, $K(z,w)$ is anti-analytic for $w$ in a neighbourhood of $\overline{\D};$\par
\smallskip
\noindent (ii) The kernel $K(z,w)$ satisfies \emph{the infinite-norm-property at the boundary}, that is,
 \begin{eqnarray}\label{infinite norm}
 \lim_{w\to \partial \D}\|K_w\|=\infty;
 \end{eqnarray}
 and,\par
 \smallskip
 \noindent (iii) $K(z,w)$ satisfies
 \emph{the uniformly boundedness condition}
 \begin{eqnarray}\label{bdd}
 \frac{|K_w(z)|}{\|K_w\|^2}\leq C_{\mathcal{H}},\quad w, z\in \D,
 \end{eqnarray}
 where $C_{\mathcal{H}}$ is a constant depending on the space.\par
 \smallskip
  There is an ordered sub-family, $\mathbb{H}_{W_\beta}({\D}), -\infty<\beta<\infty,
 $ called the Hardy-Sobolev spaces, within the family of weighted Hardy spaces (\cite{Ball1,Ball2}, also see \S 4). The index range $\beta<0$ corresponds to the weighted Bergman spaces including the standard Bergman space case for $\beta=-1$ whose $n$-best existence results are proved in \cite{QQLZ}. The $\beta=0$ case corresponds to the Hardy space \cite{WQ2020}. The $n$-best existence results for the Hardy-Sobolev spaces for $0<\beta\leq 1$ ($\beta=1$ corresponds to the Dirichlet space) are obtained as a consequence of the main result of this paper (see \S 4).  The Hardy-Sobolev spaces for $\beta>1$ do not fall into the category of the RKHSs considered in the main theorem of this paper, but we show that they are governed by the Sobolev Embedding Theorem (also see \cite{qu2019}). \\

 The writing plan of the paper is as follows. In \S 2 we discuss in detail the Gram-Schmidt orthogonalization of reproducing kernels that induces the concept multiple kernels. It is in terms of the multiple kernel concept that the $n$-best problem is precisely formulated.  In \S 3 our main result, Theorem \ref{Theorem 1}, on existence of the $n$-best approximation of a wide class of RKHSs is proved  through a number of technical lemmas in relation to orthogonal projections and analysis of the involved reproducing kernels. In \S 4, as examples of using the main result Theorem \ref{Theorem 1}, we give re-proofs of the Hardy and the weighted Bergman space results of \cite{WQ2020} and \cite{QQLZ}, and to prove, using the unified method, the $n$-best existence result for the range $-\infty<\beta\leq 1$ of the Hardy-Sobolev spaces $\mathbb{H}_{W_\beta}({\D}),$ in which the results for the range $(-\infty,1]$ are known to be equivalent to the weighted Bergman and the Hardy space cases. The results for the range $\beta\in (0,1]$ are new as applications of Theorem \ref{Theorem 1}. We include a remark in \S 4 concerning the range $\beta\in (1,\infty)$ through invoking the Sobolev Embedding theorem that provides complete understanding to the $n$-best issue for the whole range $-\infty <\beta <\infty.$ In \S 5 we extend Theorem \ref{Theorem 1} to the stochastic signal case based on the Bochner type Hilbert space setting. For the existing related studies in the stochastic signal direction we refer to \cite{QianSAFD} and \cite{QQD}. To the end of \S 5 we include a remark for the impact of this study on obtaining algorithms in finding the $n$-best solutions.

\section{Main Theorem}

Let $\mathcal{H}$ be a RKHS of holomorphic functions in $\D$ with reproducing kernel $K_w, w\in \D,$ satisfying the conditions (i),(ii),(iii) set in \S 1.  Throughout the paper $n$ is a fixed positive integer. Let $Z_k=(a_1,\cdots,a_k), 1\leq k\leq n,$ be an ordered $k$-tuple of complex numbers in $\D$ allowing multiplicity.

Denote by $l(a_k)$ the multiple of $a_k$ in the $k$-tuple $(a_1,\cdots,a_k), k\leq n.$
Denote
\begin{eqnarray}
\tilde{K}_{a_k}(z)=\left[\left(\frac{d}{d\overline{w}}\right)^{l(a_k)-1}
K_w(z)\right]_{w=a_k}.
\label{derivative}
\end{eqnarray}
We will call $\tilde{K}_{a_k}$
the {\it multiple reproducing kernel  corresponding to}  $(a_1,\cdots,a_k)$.
It is easy to show, for $f$ being in the holomorphic function space, there holds
\begin{eqnarray}\label{beautiful} \langle f, \tilde{K}_{a_k}\rangle =f^{(l(a_k)-1)}(a_k). \end{eqnarray} The consecutive derivatives of the kernel function correspond to repeating use of kernel parameters.

In the Hardy space case, for instance, $k_a(z)=\frac{1}{1-\overline{a}z},$ and
\begin{eqnarray}\label{Hardy multiple kernel}
\tilde{k}_{a_k}(z)=\left[\left(\frac{d}{d\overline{w}}\right)^{l(a_k)-1}
k_w(z)\right]_{w=a_k}=\frac{l!\overline{a}^l}{(1-\overline{a}z)^{l+1}}.
\label{derivative}
\end{eqnarray}
In general cases, let $(a_1,\cdots, a_n)$ be any $n$-tuple of complex numbers in $\D.$
Denote by
$(E_{a_1},\cdots,E_{a_1\cdots a_m})$ the Grand-Schmidt orthonormalization of $(\tilde{K}_{a_1},\cdots,\tilde{K}_{a_m}),\ m=1,\cdots,n,$ given by
\begin{eqnarray}\label{multiple}
E_{a_1\cdots a_m}(z)=\frac{\tilde{K}_{a_{m}}(z)-\sum_{l=1}^{m-1} \langle \tilde{K}_{a_{m}}, E_{a_1\cdots a_l}\rangle E_{a_1\cdots a_l} (z)}{\sqrt{\|\tilde{K}_{a_{m}}\|^2-\sum_{l=1}^{m-1}|\langle \tilde{K}_{a_{m}}, E_{a_1\cdots a_l}\rangle|^2}}.
\end{eqnarray}
We will denote the orthogonal projection of $f$ into the linear subspace $X$ by
$P_X(f).$ The projection into the orthogonal complement of $X$ is denoted $Q_X=I-P_X.$
In particular, denote by $P_{a_1\cdots a_m}$ the orthogonal projection from $\mathcal{H}$ to
${\rm span}\{\tilde{K}_{a_1},\cdots,\tilde{K}_{a_m}\},$ and by $Q_{a_1\cdots a_m}=I-P_{a_1\cdots a_m},$ the projection into the  orthogonal complement subspace of ${\rm span}\{\tilde{K}_{a_1},\cdots,\tilde{K}_{a_m}\}.$
It is recognized that $Q_{a_1a_2...a_k}$ corresponds to the Gram-Schmidt process, precisely,
 \begin{eqnarray}\label{as given} E_{a_1a_2...a_k}=
 \frac{Q_{a_1a_2...a_{k-1}}(\tilde{K}_{a_k})}{\|Q_{a_1a_2...a_{k-1}}(\tilde{K}_{a_k})\|}.\end{eqnarray}

We have been using the notation $\{B_{a_1a_2,...a_k}\}_{k=1}^\infty$ for the TM system in the Hardy space case. Now we use
$\{E_{a_1a_2...a_k}\}_{k=1}^\infty$ for the Gram-Schmidt orthonormalization of the multiple reproducing kernels $\{\tilde{K}_{a_k}\}_{k=1}^\infty$ in $\mathcal{H}$ given by (\ref{multiple}) and (\ref{as given}). In the classical Hardy space case they are essentially the same, as, in fact, $E_{a_1a_2...a_k}=c_kB_{a_1a_2...a_k},$ where $c_k$ are unimodular constants, $k=1,\cdots,$ (see Lemma \ref{lemma2} and the relevant references).  None of the $E_{a_1a_2...a_k}$ of any holomorphic Hilbert spaces other than the Hardy space seem to have such nice construction: the orthonormalization of the Szeg\"o kernel $k_{a_k}$ with respect to the span of $\{B_{a_1a_2...a_{j-1}}\}_{j=1}^{k-1},$ which is $B_{a_1a_2...a_k},$ is just the product of the added normalized Szeg\"o kernel $e_{a_k}$ with the the canonical Blaschke $\phi_{a_1a_2...a_{k-1}}.$ This extraordinary property, together with the complex unimodular property of Blaschke products on the circle, as well as the equivalent norm property restricted to the circle, offer decisive conveniences in developing the Hardy space theory in contrast with that of the non-Hardy space cases.  The $n$-best kernel approximation problem in the general context is formulated as follows. Let $f\in \mathcal{H}.$ Whether there exist, and if yes, how to find computationally $a_1,\cdots,a_n,$ all in $\D,$ such that
\begin{eqnarray}\label{hold1}
d_f=d_f(n)\triangleq\|f-P_{a_1\cdots a_n}f\|=\inf\{\|f-P_{b_1\cdots b_n}f\|\ :\ (b_1\cdots b_n) \in {\D}^n\}?
\end{eqnarray}

We will prove the following

\begin{theorem}\label{Theorem 1} Suppose that $\mathcal{H}$ is a RKHS of holomorphic functions satisfying (i), (ii), and (iii), and $n$ is a positive integer. Then there must hold one of the following two cases: (1) $f$ is a linear combination of
$\tilde{K}_{b_1},\cdots,\tilde{K}_{b_{m_1}}$ for some $m_1$-tuple $(b_1,\cdots,b_{m_1})\in {\D}^{m_1}, m_1\leq n;$ or (2) there exists an $n$-tuple $(a_1,\cdots,a_n)\in {\D}^n$ such that (\ref{hold1}) holds for a positive infimum $d_f>0,$ that is,
\begin{eqnarray}\label{holds}
d_f&=&\|f-\sum_{l=1}^n\langle f,E_{a_1\cdots a_l}\rangle E_{a_1\cdots a_l}\|\nonumber \\
&=& \inf \{\|f-\sum_{l=1}^n\langle f,E_{b_1\cdots b_l}\rangle E_{b_1 \cdots b_l}\|\ :\ (b_1,\cdots,b_n)\in {\D}^n\}.
\end{eqnarray}
 \end{theorem}

 We note that $(a_1,\cdots,a_n)\in {\D}^n$ gives rise to
  equality (\ref{holds})
 if and only if
 \begin{eqnarray}\label{sup} \sum_{l=1}^n|\langle f,E_{a_1\cdots a_l}\rangle|^2=\sup \{\sum_{l=1}^n|\langle f,E_{b_1\cdots b_l}\rangle|^2, (b_1,\cdots,b_n)\in {\D}^n\}. \end{eqnarray}

\section{Proof of the Theorem}

We will use the notation $f_{a_1a_2...a_k}=Q_{a_1a_2...a_k}f,$ where $a_j, j=1,\cdots,k,$ are allowed to repeat.
\begin{lemma}\label{multiplicity}
If $a_j$ is among $a_1,\cdots,a_k,$ then
$f_{a_1a_2...a_k}(a_j)=0,$
including the multiplicity.
\end{lemma}

\noindent{\bf Proof.} The proof is straightforward if $a_1,\cdots,a_k$ are all different, that is $l(a_j)=1, j=1,\cdots,k.$ In the case,
owing to the self-adjoint property of the projection operators and the orthogonality gained from G-S process, for $a_j$ being among $a_1,\cdots,a_k,$
\begin{eqnarray*}
 f_{a_1a_2...a_k}(a_j)&=&\langle Q_{a_1a_2...a_k}f, K_{a_j}\rangle\\
 &=&\langle f, Q_{a_1a_2...a_k}K_{a_j}\rangle\\
 &=&\langle f, (I-P_{a_1a_2...a_k})K_{a_j}\rangle\\
 &=& \langle f, 0\rangle\\
 &=& 0.\end{eqnarray*}
 Let $a_j$ have multiplicity $l(a_j)>1,$ and $a_{s_1}=a_{s_2}=\cdots=a_{s_{l(a_j)}}=a_j.$ For $m=1,\cdots,l(a_j),$ in view of (\ref{beautiful}),
 \begin{eqnarray*}
 \left(\frac{d}{dz}\right)^{m-1}[f_{a_1a_2...a_k}](a_j)&=&\langle \left(\frac{d}{dz}\right)^{m-1}Q_{a_1a_2...a_k}f, K_{a_j}\rangle\\
 &=&\langle f, Q_{a_1a_2...a_k}\left(\frac{d}{d\overline{w}}\right)^{m-1}
K_w(z)|_{w=a_j}\rangle\\
 &=&\langle f, (I-P_{a_1a_2...a_k})\tilde{K}_{a_j}\rangle\\
 &=& 0.\end{eqnarray*}
 So, $f_{a_1a_2...a_k}$ has $l(a_j)$-multiple zero at $a_j.$

 $\hfill\square$\\
\def\e{\rm e}

The new crucial concepts of this methodology include {\it higher order generalized backward shift operator}
\[  \frac{Q_{a_1\cdots{a_{m-1}}}(f)(z)}{\phi_{a_1\cdots{a_{m-1}}}(z)}\]
and its factorization (see Lemma \ref{lemma2} below).
The order-1 generalized backward shift operator
$Q_{a}/[\frac{z-a}{1-\overline{a}z}]$ applied to a test function $f$  gives rise to what we call reduced remainder (\cite{QWa}) playing a central role in the formulation of AFD.

\begin{lemma}\label{lemma2}
Let $a_1,\cdots,a_{m}$ be complex numbers in $\D$ allowing multiplicity. Then\\

(1) In the general $\mathcal{H}$ space setting the Gram-Schmidt orthonormalization
\begin{eqnarray}\label{tildeBla} E_{a_1\cdots a_{m}}(z)= \frac{Q_{a_1\cdots{a_{m-1}}}({K}_{a_m})(z)}{\|Q_{a_1\cdots{a_{m-1}}}({K}_{a_m})\|}
&=&\frac{{\tilde{K}}_{a_m}(z)-\sum_{l=1}^{m-1}\langle {\tilde{K}}_{a_m},E_{a_1\cdots a_{l}}\rangle E_{a_1\cdots a_{l}}(z)}{\|{\tilde{K}}_{a_m}-\sum_{l=1}^{m-1}\langle {\tilde{K}}_{a_m},E_{a_1\cdots a_{l}}\rangle E_{a_1\cdots a_{l}}\|}.\end{eqnarray}

(2) In the Hardy space case, the
 last function (\ref{tildeBla}) is equal to ${\e}^{ic}e_{a_m}(z)\phi_{a_1\cdots a_{m-1}}(z)
={\e}^{ic}B_{a_1\cdots a_{m}}(z),$
where
\[ {\e}^{ic}=\frac{\overline{\phi}_{a_1\cdots{a_{m-1}}}(a_m)}
{|\overline{\phi}_{a_1\cdots{a_{m-1}}}(a_m)|},\]
 and $e_{a_m},\phi_{a_1\cdots{a_{m-1}}}$ and $B_{a_1\cdots a_{m}}$ are respectively defined in (\ref{e}), (\ref{phi}) and (\ref{B}). Moreover,
\begin{eqnarray}
 \langle f,E_{a_1\cdots a_m}\rangle E_{a_1\cdots a_m}=\langle f,B_{a_1\cdots a_m}\rangle B_{a_1\cdots a_m},\end{eqnarray}
 and
 \begin{eqnarray}
 \langle f,E_{a_1\cdots a_m}\rangle=\frac{Q_{a_1\cdots{a_{m-1}}}(f)(a_m)}{\phi_{a_1\cdots{a_{m-1}}}(a_m)}\sqrt{1-|a_m|^2}.
 \end{eqnarray}

 (3) In the general $\mathcal{H}$ space setting, the higher order generalized shift operators (the reduced remainders) may be factorized, that is

 \begin{eqnarray}\label{iteration}
 \frac{Q_{a_1\cdots{a_{m-1}}}(f)(z)}{\phi_{a_1\cdots{a_{m-1}}}(z)}=
 \left(\frac{Q_{a_{m-1}}}{\phi_{a_{m-1}}}
 \circ\cdots\circ\frac{Q_{a_1}}{\phi_{a_1}}\right)(f)(z).
 \end{eqnarray}
\end{lemma}

\noindent{\bf Proof.} For proofs of (1) and (2) we refer to \cite{QianSAFD}. We now prove (3).
For $k>1,$ denote by $g_k$ the \emph{ $k$-reduced remainder} (\cite{QWa})
\[ g_{k+1}(z)=\frac{g_{k}(z)-\langle g_{k},E_{a_{k}}\rangle E_{a_{k}}(z)}{\phi_{a_{k}}(z)}=\left(\frac{Q_{a_k}}{\phi_{a_k}}\right)(g_k)(z),\]
where $g_1=f.$
Inductively there holds
\[g_{k+1}(z)=\left(\frac{Q_{a_k}}{\phi_{a_k}}\circ \frac{Q_{a_{k-1}}}{\phi_{a_{k-1}}}\right)(g_{k-1})(z)=\left(
\frac{Q_{a_k}}{\phi_{a_k}}\circ \frac{Q_{a_{k-1}}}{\phi_{a_{k-1}}}\circ
\cdots \frac{Q_{a_1}}{\phi_{a_1}}\right)(f)(z).\]

On the other hand, the AFD formulation given in \cite{QWa} implies
\[Q_{a_1\cdots a_{k}}f=g_{k+1}\prod_{j=1}^k \phi_{a_j}.\]
We thus have
\[g_{k+1}=\frac{Q_{a_1\cdots a_{k}}f}{\prod_{j=1}^k \phi_{a_j}}=\left(
\frac{Q_{a_k}}{\phi_{a_k}}\circ \frac{Q_{a_{k-1}}}{\phi_{a_{k-1}}}\circ
\cdots \frac{Q_{a_1}}{\phi_{a_1}}\right)(f).\]
$\hfill\square$\\

\begin{remark}
It is recognized that the operator $\frac{Q_a}{\phi_a}$ is the generalized backward shift operator defined in \cite{QWa}. Repeating use of the operator yields the reduced remainders $\frac{Q_{a_1\cdots a_{k}}f}{\prod_{j=1}^k \phi_{a_j}}.$ The following lemma shows that the reduced remainders of a function bounded by $M$ are still bounded with explicit bounds in terms of $M$ and the parameters $a_1,\cdots,a_k$ involved. The result plays a crucial role in the proof of the main result of the paper.
\end{remark}

 \begin{lemma}\label{Mbdd}
 Let $f$ be an analytic function in an open neighbourhood of $\overline{\D}$ and
 \[ |f(z)|\leq M\]
 for some $M>0$ on $\overline{\D}.$
 Then for any sequence $a_1,a_2,...,a_k\in \D$ the reduced remainder functions
 \[ \frac{f_{a_1a_2...a_k}(z)}{\phi_{a_1a_2...a_k}(z)}\]
 are analytic in an
  open neighbourhood of $\overline{\D}$ with the bounds over $\overline{\D}:$
 \begin{eqnarray} \label{repeating} \left|\frac{f_{a_1a_2...a_k}(z)}{\phi_{a_1a_2...a_k}(z)}\right|\leq M(1+C_{\mathcal{H}})^k,\quad z\in \overline{\D},\end{eqnarray}
 where $C_{\mathcal{H}}$ is the constant in (\ref{bdd}).
 \end{lemma}

\noindent {\bf Proof.} For $a_1\in \D$ and $z\in \partial \D,$
\begin{eqnarray*} |f_{a_1}(z)|&=&|f(z)-\langle f,E_{a_1}\rangle E_{a_1} (z)|\\
&\leq& |f(z)|+|f(a_1)|\frac{|K_{a_1}(z)|}{\|K_{a_1}\|^2}\\
&\leq& |f(z)|+|f(a_1)|C_{\mathcal{H}}\\
&\leq& M(1+C_{\mathcal{H}}).
\end{eqnarray*}
Since the zero of $\phi_{a_1}$ is a zero of $Q_{a_1}f,$ $\frac{Q_{a_1}f}{\phi_{a_1}}$ is a holomorphic function in an open neighbourhood of $\overline{\D}.$  The maximum modulus principle over $\overline{\D}$ gives
\[ \left|\frac{Q_{a_1}f(w )}{\phi_{a_1}(w )}\right|\leq \max\{|f_{a_1}(z)|\ :\ z\in {\partial \D}\}\leq M(1+C_{\mathcal{H}})\]
for all $w\in \overline{\D}.$
By invoking the result of (iii), Lemma \ref{lemma2}, and repeating $k$ times the above estimation for $\frac{Q_{a_1}f}{\phi_{a_1}},$ we obtain the bounds claimed by the lemma.
 $\hfill\square$\\

We first prove Theorem 2.1 for $n=1.$ The goal of $1$-best approximation amounts to finding $a_1\in \D$ such that \[|\langle f,E_{a_1}\rangle|=\frac{|f(a_1)|}{\|K_{a_1}\|}\]
attains its global maximum over all possible choices of the parameter in $\D.$ This is what we called  Maximal Selection Principle (MSP) in the previous related studies (\cite{QWa,qu2018,qu2019}). The proof is divided into two steps:\\

\noindent {\bf Step 1}: We show that Boundary Vanishing Condition (BVC) holds in the case, that is,
\begin{eqnarray}\label{BVC} \lim_{|a|\to 1-}|\langle f,E_{a}\rangle |=0.\end{eqnarray}

BVC is proved here by using a density argument: We note that the given function $f\in \mathcal{H}$ may be approximated within an error $\epsilon>0$ by a bounded holomorphic function $g$ as a linear combination of some parameterized reproducing kernels.
Therefore, by the Cauchy-Schwarz inequality, we have
 \[ |\langle f,E_a\rangle|\leq |\langle f-g,E_a\rangle|+ |\langle g,E_a\rangle|\leq
 \|f-g\|+\frac{|g(a)|}{\|K_a\|}\leq \epsilon +\frac{|g(a)|}{\|K_a\|}.\]
As a consequence of (\ref{infinite norm}) and boundedness of $g$, the BVC  (\ref{BVC}) is concluded.\par
\smallskip
\noindent {\bf Step 2}: If $f$ is not identical with the zero function there exists $b\in \D$ such that $|\langle f,E_b\rangle|>0.$ Denote $|\langle f,E_b\rangle|=\delta.$ The argument in Step 1 amounts that there exists $0<r_1<1$ such that $|a|>r_1$ implies $|\langle f,E_a\rangle|<\delta/2.$ Therefore,
 \begin{eqnarray}\label{show} \max \{|\langle f,E_a\rangle|\ :\ |a|\leq r_1\}=\sup \{|\langle f,E_a\rangle|\ :\ a\in \D\}.\end{eqnarray}
 By invoking the theorem of global maximum attainability of a continuous function on a compact set, (\ref{show}) shows that the global maximum of $|\langle f,E_a\rangle|$ is attainable inside $\D$.\\

%

Now we turn directly to the general $n>1$ case. To begin with the proof we assume that $f$ itself is not expressible as a linear combination of $m_1$ reproducing kernels for $m_1<n.$ Based on the definition of supreme, one can find a sequence of $n$-tuples with mutually distinct and non-zero components (these can always be done owing to continuity of the inner product, and can at least simplify the notation),
$(a_1^{(l)},\cdots,a_n^{(l)}), l=1,2,\cdots,$ that corresponds to a sequence of $n$-tuples of reproducing kernels $(K_{a_1^{(l)}},\cdots,K_{a_n^{(l)}}),$ such that the norms of the projections $P_{a_1^{(l)}\cdots a_n^{(l)}}(f)$ tends to the supreme
 (\ref{sup}). Since $(a_1^{(l)},\cdots,a_n^{(l)})\in\overline{\D}^n,$ through a Bolzano-Weierstrass compact argument, we may assume, without loss of generality, that the sequence of the $n$-tuples $(a_1^{(l)},\cdots,a_n
 ^{(l)})$ itself converges to $(a_1,\cdots,a_n)\in \overline{\D}^n.$ If we have $a_1,\cdots,a_n$ all in ${\D},$ then we are done due to continuity of the inner product, although may involve multiple kernels when multiplicities occur. This gives rise to the case (1) of the Theorem \ref{Theorem 1} for $m_1=n$ when $d_f=0;$ and the case (2) when $d_f>0.$ \par
 \smallskip
 Now we show that if not all the limiting points $a_1,\cdots,a_n$ are located within ${\D},$ then we have the case (1) for some $m_1<n,$ being contrary with our priori assumption. Assume that at least one of $a_1,\cdots,a_n$ are on the boundary $\partial \D.$ Since the projections $P_{a_1^{(l)}\cdots a_n^{(l)}}(f)$ and $Q_{a_1^{(l)}\cdots a_n^{(l)}}(f)$ are irrelevant with the order, by re-ordering, when necessary, we may assume without loss of generality that $a_1,\cdots,a_{m_1}$ are in $\D,$ and $a_{m_1+1},\cdots, a_n$ are on $\partial \D,$ where $m_1<n,$ and in particular, $\lim_{l\to\infty} |a_n^{(l)}|=1.$
 In the case we will show
 \begin{eqnarray}\label{rest} \lim_{l\to \infty} |\langle f,E_{a^{(l)}_1a^{(l)}_2...a^{(l)}_n}\rangle|=0,\end{eqnarray} regardless
 the locations of $a^{(l)}_k, k=1,\cdots,n-1$ and $l=1,2,\cdots$. If (\ref{rest}) can be proved,
 by repeating the same argument $n-m_1$ times we result in that the latter $n-m_1$ terms of the $l$-sequence of the $n$-tuples all have no contribution. We claim that $d_f>0$ cannot hold. If it were $d_f>0,$ then it could be further reduced, contrary with $d_f$ being infimum. But, $d_f$ cannot be zero either, for in such case we got that $f$ is a linear combination of $m_1<n$ multiple reproducing kernels, contrary with our assumption. Thus, all that remain to be proved is (\ref{rest}). By using the same density argument as we prove the case $n=1$ we may assume that $f$ itself is an analytic function in a neighbourhood of the closed unit disc
$\overline{\D}$ with a bound $M.$

 For any zero set $Z$ possibly with multiplicities we use the general notation  $K_{Z}(z,a)$ for the reproducing kernel at $a$ of the zero space $\mathcal{H}_{Z},$ where
\[\mathcal{H}_{Z}=\{f\in \mathcal{H}\ :\ f \ {\rm vanishes\ at\ points \ in}\ Z\ {\rm including \ multiplicities}\}.\]
The space $\mathcal{H}_{Z}$ uses the same inner product as $\mathcal{H}.$
We denote by $\mathcal{H}_{\phi_Z}$ the Hilbert space
\[ \mathcal{H}_{\phi_Z}=\{ f:{\D}\to {\C}\ :\ f\ {\rm is\ analytic},\ \|f\phi_Z\|_{\mathcal{H}}<\infty\},\]
where $\phi_Z$ is the canonical Blaschke product generated by the elements of $Z$ including multiplicities.
The inner product of $\mathcal{H}_{\phi_Z}$ is denoted as $\langle \cdot,\cdot\rangle_{\mathcal{H}_{\phi_Z}}.$ The reproducing kernel of $\mathcal{H}_{\phi_Z}$ is denoted $K_{{\phi_Z}}.$
 In this paper we only need to treat zero sets $Z$ with finite points. Note that $\mathcal{H}\subset \mathcal{H}_{\phi_Z},$ and $\|f\|_{\mathcal{H}_{\phi_Z}}\leq \|f\|_{\mathcal{H}}.$\par
\smallskip
 The next two lemmas follow similar idea in \cite{zhubook1,Zhu,zhubook2,duren2004bergman,QQLZ}.

 \begin{lemma}\label{zero space kernel} For any finite zero set $Z,$ by denoting
 $K_Z(z,w)$ the reproducing kernel of the zero space $\mathcal{H}_{Z},$
  there holds
 \begin{eqnarray}\label{have} K_Z(z,w)=\phi_Z(z)K_{\phi_Z}(z,w)\overline{\phi_Z(w)}\end{eqnarray}
 and
 \begin{eqnarray}\label{simply} \|K_Z(\cdot,w)\|_{\mathcal{H}}\leq\|K_{\phi_Z}(\cdot,w)\|_{\mathcal{H}_{\phi_Z}},\end{eqnarray}
 where $\phi_Z$ is the canonical Blaschke product defined by $Z,$ ${\mathcal{H}_{\phi_Z}}$ is the $|\phi_{Z}|^2$-weighted $\mathcal{H}$ space,
 $K_{\phi_Z}$ is its reproducing kernel. As a consequence, the normalized reproducing kernel is
 \begin{eqnarray}\label{last} \frac{K_Z(z,w)}{\|K_Z(\cdot,w)\|}=\frac{\overline{\phi_Z(w)}}{|\phi_Z(w)|}
 \frac{\phi_Z(z)K_{\phi_Z}(z,w)}{K^{1/2}_{\phi_Z}(w,w)}.\end{eqnarray}
\end{lemma}

\noindent{\bf Proof.}
We note that $K_Z(z,w)$ has zero set $Z$ for the variable $z$ when $w$ is fixed in a neighbourhood of $\overline{\D};$ and zero set $Z$ for the variable $w$ when $z$ is fixed in a neighbourhood of $\overline{\D}.$   Therefore, $\phi_Z^{-1}(z)K_Z(z,w)\overline{\phi}_Z^{-1}(w)$ is analytic for $z$ in a neighbourhood of $\overline{\D},$ and anti-analytic for $w$ in a neighbourhood of $\overline{\D}.$ Let $f\in \mathcal{H}_{\phi_Z}.$ In the case $f\phi_Z\in \mathcal{H}_Z.$ We have
\begin{eqnarray*}
\langle f,\phi^{-1}_{Z}(\cdot)K_{{Z}}(\cdot,w)\overline{\phi}^{-1}_{Z}(w)
\rangle_{\mathcal{H}_{\phi_Z}}&=& \phi^{-1}_{Z}(w)\langle f\phi_{Z},K_{{Z}}(\cdot,w)\rangle_{\mathcal{H}_{Z}}\\
&=&\phi^{-1}_{Z}(w)\langle f\phi_Z,K_{{Z}}(\cdot,w)\rangle_{\mathcal{H}_{Z}}\\
&=&\phi^{-1}_{Z}(w)f(w)\phi_Z(w)\\
&=&f(w).
\end{eqnarray*}
Therefore, $\mathcal{H}_{\phi_Z}$ is a RKHS. Due to uniqueness of reproducing kernel,
 its kernel $K_{\phi_Z}(z,w)$ satisfies the relation (\ref{have}). To prove (\ref{simply}) we have
\begin{eqnarray*}
\|K_Z(\cdot,w)\|^2_{\mathcal{H}}&=&K_Z(w,w)\\
&=& \phi_{Z}(w)K_{\phi_{Z}}(w,w)\overline{\phi_{Z}(w)}\\
&=& K_{\phi_{Z}}(w,w)|\phi_{Z}(w)|^2\\
&\leq& K_{\phi_{Z}}(w,w)\\
&=& \|K_{\phi_{Z}}(\cdot,w)\|^2_{\mathcal{H}_{\phi_{Z}}}.
\end{eqnarray*}
The relation (\ref{last}) is just by dividing $K_Z(z,w)$ with $\|K_Z(\cdot,w)\|=
\sqrt{K_Z(w,w)}$ and invoking (\ref{have}).\\

 $\hfill\square$

\begin{lemma}\label{zero space norm}
For any Hilbert space $\mathcal{H}$ with reproducing kernel $K$ and any $a\in\D$ there hold
\begin{eqnarray}\label{H}
K(a,a)=\sup\{ |f(a)|^2\ :\ f\in \mathcal{H}, \|f\|\leq 1\}
\end{eqnarray}
and
\begin{eqnarray}\label{normal}
K(a,a)\leq K_{\phi_Z}(a,a).
\end{eqnarray}
\end{lemma}
  \noindent{\bf Proof.} Recall that $E_a(z)=K(z,a)/\|K_a\|.$ On one hand, $\|E_a\|=1.$ On the other hand, for any $f$ satisfying $\|f\|=1,$ using the Cauchy-Schwarz inequality,
  \[ |f(a)|^2=|\langle f,K_a\rangle|^2\leq \|K_a\|^2=K(a,a).\]
  So, $E_a$ is a solution for the extremal problem.
  Using this argument also to $\mathcal{H}_{\phi_{Z}}$ and $K_{\phi_{Z}}$, we obtain
  \begin{eqnarray*}
  K(a,a)&=&\sup\{|f(a)|^2\ :\ f\in \mathcal{H}, \|f\|_{\mathcal{H}}\leq 1\}\\
  &\leq&\sup\{|f(a)|^2\ :\ f\in \mathcal{H}_{\phi_Z}, \|f\|_{\mathcal{H}_{\phi_Z}}\leq 1\}\\
  &=&K_{{\phi_{Z}}}(a,a),
  \end{eqnarray*}
  as desired.
   $\hfill\square$\\

  Now we proceed with the main technical step of the proof.
Denote by $Z_{n-1}^{(l)}$ the $l$-level zero set $(a^{(l)}_1,a^{(l)}_2,...,a^{(l)}_{n-1}).$
With the above preparations we have
\begin{eqnarray*}\label{last term}
  & & \langle f,E_{a^{(l)}_1 a^{(l)}_2...a^{(l)}_n}\rangle_{\mathcal{H}} \\
  &=& \langle f_{a^{(l)}_1 a^{(l)}_2...a^{(l)}_{n-1}},
  E_{a^{(l)}_1 a^{(l)}_2...a^{(l)}_n}\rangle_{\mathcal{H}} \qquad \left(E_{a^{(l)}_1 a^{(l)}_2...a^{(l)}_n}=\frac{Q_{a^{(l)}_1 a^{(l)}_2...a^{(l)}_{n-1}}(K_{a^{(l)}_n})}{\|Q_{a^{(l)}_1 a^{(l)}_2...a^{(l)}_{n-1}}(K_{a^{(l)}_n}\|}=\frac{Q^2_{a^{(l)}_1 a^{(l)}_2...a^{(l)}_{n-1}}(K_{a^{(l)}_n})}{\|Q_{a^{(l)}_1 a^{(l)}_2...a^{(l)}_{n-1}}(K_{a^{(l)}_n}\|}\right)\nonumber \\
&=& \left\langle f_{a^{(l)}_1 a^{(l)}_2...a^{(l)}_{n-1}}, \frac{K_{Z_{n-1}^{(l)}}(\cdot,a^{(l)}_n)}{\|K_{Z_{n-1}^{(l)}}(\cdot,a^{(l)}_n)\|}\right\rangle_{\mathcal{H}}
\qquad \left(Q_{a^{(l)}_1 a^{(l)}_2...a^{(l)}_{n-1}}(K_{a^{(l)}_n})=K_{Z_{n-1}^{(l)}}(\cdot,a^{(l)}_n)\right)\\
  &=& \frac{\overline{\phi_{Z^{(l)}_{n-1}}(a^{(l)}_n)}}{|\phi_{Z^{(l)}_{n-1}}(a^{(l)}_n)|}\langle f_{a^{(l)}_1 a^{(l)}_2...a^{(l)}_{n-1}}, K_{\phi_{Z_{n-1}^{(l)}}} (a^{(l)}_n,a^{(l)}_n)^{-1/2}\phi_{Z^{(l)}_{n-1}}K_{\phi_{Z_{n-1}^{(l)}}} (\cdot ,a^{(l)}_n) \rangle_{\mathcal{H}} \qquad \left({\rm Lemma}\ \ref{zero space kernel}\right)\\
  &=&\frac{\overline{\phi_{Z^{(l)}_{n-1}}(a^{(l)}_n)}}{|\phi_{Z^{(l)}_{n-1}}(a^{(l)}_n)|}
  \left\langle \frac{f_{a^{(l)}_1 a^{(l)}_2...a^{(l)}_{n-1}}}{\phi_{Z^{(l)}_{n-1}}}, |\phi_{Z^{(l)}_{n-1}}|^2 K_{\phi_{Z_{n-1}^{(l)}}} (\cdot ,a^{(l)}_n)
  \right\rangle_{\mathcal{H}}\frac{1}{\sqrt{K_{\phi_{Z_{n-1}^{(l)}}} (a^{(l)}_n,a^{(l)}_n)}}\qquad \left({\rm Lemma}\ \ref{zero space kernel}\right)\\
  &=& \frac{\overline{\phi_{Z^{(l)}_{n-1}}(a^{(l)}_n)}}{|\phi_{Z^{(l)}_{n-1}}(a^{(l)}_n)|}
  \left\langle \frac{f_{a^{(l)}_1 a^{(l)}_2...a^{(l)}_{n-1}}}{\phi_{Z^{(l)}_{n-1}}}, K_{\phi_{Z_{n-1}^{(l)}}} (\cdot ,a^{(l)}_n)
  \right\rangle_{\mathcal{H}_{\phi_{Z^{(l)}_{n-1}}}}
  \frac{1}{\sqrt{K_{\phi_{Z_{n-1}^{(l)}}} (a^{(l)}_n,a^{(l)}_n)}}\\
&=&\frac{\overline{\phi_{Z^{(l)}_{n-1}}(a^{(l)}_n)}}{|\phi_{Z^{(l)}_{n-1}}(a^{(l)}_n)|}
\frac{f_{a^{(l)}_1 a^{(l)}_2...a^{(l)}_{n-1}}(a^{(l)}_n)}
{\phi_{Z^{(l)}_{n-1}}(a^{(l)}_n)}
\frac{1}{\sqrt{K_{\phi_{Z_{n-1}^{(l)}}} (a^{(l)}_n,a^{(l)}_n)}}\quad \left({\rm Lemma}\ \ref{zero space kernel}\ {\rm and\ reproducing \ kernel\ property} \right).
\end{eqnarray*}
To conclude the theorem it is sufficient to show that the above quantity tends to zero along with ${\D}\ni a^{(l)}_n\to a_n\in \partial \D$ uniformly in $a^{(l)}_1,\cdots,a^{(l)}_{n-1} \in \D$ for $l=1,2\cdots$ It then suffices to prove \\

$1^o.$
\[ \frac{f_{a^{(l)}_1 a^{(l)}_2...a^{(l)}_{n-1}}(a^{(l)}_n)}
{\phi_{Z_{n-1}^{(l)}} (a^{(l)}_n)}\]
is bounded uniformly in $a_1^{(l)},\cdots,a_{n-1}^{(l)}$ and $a_{n}^{(l)}, l=1,2,\cdots;$ and\\

$2^o.$
\[ \lim_{l\to \infty}K_{\phi_{Z_{n-1}^{(l)}}} (a^{(l)}_n,a^{(l)}_n)=\infty\]
uniformly in $a_1^{(l)},\cdots,a_{n-1}^{(l)}, l=1,2,\cdots$\par
\smallskip
Now we show assertion $1^o.$ First by Lemma \ref{multiplicity} the function
\[ g_{Z^{(l)}_{n-1}}(z)=\frac{f_{a^{(l)}_1\cdots a^{(l)}_{n-1}}(z)}{\phi_{Z^{(l)}_{n-1}}(z)}\]
is analytic in a neighbourhood of $\overline{\D}.$ By invoking the maximum modulus principle for one complex variable in $\overline{\D},$ Lemma \ref{Mbdd}, as well as the fact that all finite Blaschke products are of modulus $1$ on the boundary $\partial{\D},$ we have
\begin{eqnarray*}
\max\{|g_{Z^{(l)}_{n-1}}(z)|\ :\ z\in \overline{\D}\}&=&\max\{|g_{Z^{(l)}_{n-1}}(\zeta)|\ :\ \zeta\in \partial{\D}\}\\
&=&\max\{\left| f_{a^{(l)}_1 \cdots a^{(l)}_{n-1}}({\rm e}^{it}) \right|\ :\ t\in \partial{\D}\}\\
&\leq& M(1+C_{\mathcal{H}})^{n-1},
\end{eqnarray*}
concluding the uniform boundedness claimed of $1^o.$  The assertion $2^o$ is a consequence of the condition (\ref{infinite norm}) and Lemma \ref{zero space norm}.

\section{Applications}

\subsection{The classical Hardy space}
By taking $\mathcal{H}=\mathbb{H}^2({\D}),$ we are with the inner product (\ref{inner}), and
the reproducing kernel $K(z,w)=k_w(z)=\frac{1}{1-\overline{w}z}.$
Since
\[ K(a,a)=\frac{1}{1-|a|^2}\to \infty \ {\rm as}\ |a|\to 1,\ {\rm and}\ \frac{|K_a(z)|}{K(a,a)}=\frac{1-|a|^2}{|1-\overline{a}z|}\leq 2,\]
the conditions (\ref{infinite norm}) and (\ref{bdd}) are satisfied. We hence have existence of the $n$-best approximation.

\subsection{The Bergman spaces}

Let $\mathcal{H}$ be the weighted Bergman spaces with the definition and notation
\begin{eqnarray*}
\mmA({\D})=\{f: {\D}\to {\bf C}\ | \ f\ {\rm is \ holomorphic \ in\ \D}, {\rm and}\ \|f\|_{\mmA({\D})}^2=\int_{\D} |f(z)|^2dA_\alpha<\infty\},
\end{eqnarray*}
where $\alpha\in (-1,\infty), dA_\alpha=(1+\alpha)(1-|z|^2)^{\alpha}dA(z),$ and $dA=\frac{dxdy}{\pi}, z=x+iy,$ is the normalized area measure of the unit disc. The inner product of $\mmA({\D})$ is defined as
$$
\langle f,g\rangle_{\mmA(\D)}=\int_{\D} f(z)\overline{g(z)}dA_\alpha.
$$
In the sequel we sometimes write $\| \cdot \|_{\mmA({\D})}$ and $\langle \cdot,\cdot\rangle_{\mmA({\D})}$ briefly as $\|\cdot \|$ and $\langle \cdot,\cdot\rangle,$ and
$\mmA({\D})$ as $\mmA.$

$\mmA$ is a RKHS with reproducing kernel
$$ k_a^\alpha(z)=\frac{1}{(1-\overline{a}z)^{2+\alpha}}.$$
By invoking the reproducing kernel property we have
\begin{eqnarray}\label{recall} \| k_a^\alpha\|^2= k_a^\alpha(a)=\frac{1}{(1-|a|^2)^{2+\alpha}}.\end{eqnarray}
This shows that the condition (\ref{infinite norm}) holds for all $\mmA.$
A simple computation gives
$$ \frac{k_a^\alpha(z)}{{k_a^\alpha(a)}}=\frac{{(1-|a|^2)^{2+\alpha}}}{(1-\overline{a}z)^{2+\alpha}}
\leq 2^{2+\alpha}.$$

Hence, the reproducing kernel satisfies the condition (\ref{bdd}). Therefore, an $n$-best approximation exists in all the weighted Bergman spaces. This is a re-proof of the main result of \cite{QQLZ}.

\subsection{The weighted Hardy spaces}

Let $W(k)$ be a sequence of non-negative numbers satisfying $\lim_{k\to\infty}W(k)^{\frac{1}{k}}\ge 1$ (\cite{Mac}). Denote by $\mathbb{H}_W({\D})$ the
$W$-weighted Hardy $\mathbb{H}^2$-space defined by
\[\mathbb{H}_W({\D})=\{f\ :\ {\D}\to {\C}\ :\ f(z)=\sum_{k=0}^\infty c_kz^k, \ z\in {\D},
\|f\|_{\mathbb{H}_W}=\sum_{k=0}^\infty W(k)|c_k|^2<\infty\}.\]
We will be considering an ordered sequence of $W$-weighted Hardy spaces defined by the weights $W_\beta(k)=(1+k)^\beta, -\infty <\beta <\infty.$ This class of function spaces is a generalization of the Hardy and the weighted Bergman spaces. In fact, $\mathbb{H}_{W_0}({\D})=\mathbb{H}^2({\D}),$ and
$\mathbb{H}_{W_\beta}({\D})=\mmA ({\D}), \alpha =-\beta-1, \beta<0 \ (\alpha>-1),$ $\mathbb{H}_{W_{-1}}({\D})$ is the standard Bergman,
and $\mathbb{H}_{W_1}({\D})$ is the Dirichlet space in $\D.$ The spaces $\mathbb{H}_{W_\beta}({\D})$ are, as a matter of fact, equivalent with the Hardy-Sobolev spaces $W^{\frac{\beta}{2},2}.$ From the last two subsections we know that the spaces $\mathbb{H}_{W_\beta}({\D}), \beta\leq 0, $ have $n$-best approximation. We now extend the result to $0<\beta\leq 1.$

The inner product of $\mathbb{H}_{W_\beta}({\D}), -\infty<\beta<\infty,$ is
\[ \langle f,g\rangle=\sum_{k=0}^\infty (k+1)^\beta c_k\overline{d}_k,\]
where $c_k$ and $d_k$ are, respectively, the coefficients of the Taylor expansions of $f$ and $g.$ From this it can be directly verified that the reproducing kernel of $\mathbb{H}_{W_\beta}({\D})$ is
\[ k^\beta_a(z)=\sum_{k=0}^\infty \frac{(z\overline{a})^k}{(k+1)^\beta}.\]
The function \begin{eqnarray}\label{owing to} k^\beta_a(a)=\sum_{k=0}^\infty \frac{|a|^{2k}}{(k+1)^\beta}\end{eqnarray}
is an increasing function in $|a|,$ and for any large $N,$
\[ \underline{\lim}_{|a|\to 1-}k^\beta_a(a)\ge \lim_{|a|\to 1-}\sum_{k=0}^N \frac{|a|^{2k}}{(k+1)^\beta}=\sum_{k=0}^N \frac{1}{(k+1)^\beta}.\]
Therefore, for all $\beta\leq 1,$
\[ \lim_{|a|\to 1-}k^\beta_a(a)=\infty,\]
verifying (\ref{infinite norm}). Next we show that the weighted Hardy spaces kernels $k^\beta_a$ satisfy the condition (\ref{bdd}). This requires to prove that the function
$\frac{|k^\beta_a(z)|}{k^\beta_a(a)}$ is uniformly bounded in $a, z\in {\D}.$
The following estimation uses the well know technique for summing up series of positive decreasing entries:
If $f$ is a positive decreasing function integrable over $(0,\infty),$ then
 \[ \int_1^\infty f(t)dt\leq \sum_{k=1}^\infty f(k)\leq \int_0^\infty f(t)dt.\]
 The estimation amounts to numerically comparing some elementary integrals. Denote by $|a|=r<1.$ Then\\
\begin{eqnarray*}
\frac{|k^\beta_a(z)|}{k^\beta_a(a)}&=&\frac{|\sum_{k=0}^\infty \frac{(\overline{a}z)^k}{(1+k)^\beta}|}{\sum_{k=0}^\infty \frac{|a|^{2k}}{(1+k)^\beta}}\\
&\leq&\frac{\sum_{k=0}^\infty \frac{r^k}{(1+k)^\beta}}{\sum_{k=0}^\infty \frac{r^{2k}}{(1+k)^\beta}}\\
&\leq&\frac{\int_0^\infty \frac{r^x}{(1+x)^\beta}dx}{\int_1^\infty \frac{r^{2x}}{(1+x)^\beta}dx}\\
&=&\frac{\int_0^\infty \frac{r^x}{(1+x)^\beta}dx}{2^{\beta-1}\int_2^\infty \frac{r^{x}}{(2+x)^\beta}dx}\qquad ({\rm change \ of\ variable})\\
&=&\frac{\int_0^\infty \frac{r^x}{(1+x)^\beta}dx}{2^{\beta-1}\left(\int_0^\infty \frac{r^{x}}{(2+x)^\beta}dx-\int_0^2 \frac{r^{x}}{(2+x)^\beta}dx\right)}\qquad (\frac{\infty}{\infty}\ {\rm type \ when\ }\ r\to 1-)\\
&\leq&\frac{\int_0^\infty \frac{r^x}{(1+x)^\beta}dx}{2^{\beta-2}\int_0^\infty \frac{r^{x}}{(1+x)^\beta}\frac{(1+x)^\beta}{(2+x)^\beta}dx}\qquad ({\rm if}\ r\ge \
{\rm some}\ r_0\in (0,1))\\
&\leq& 4.  \qquad ({\rm mean-value\ theorem \ of\ integration})
\end{eqnarray*}
For $r\leq r_0$ the estimated quantity also has a uniform bound. Hence, for $\beta\leq 1,$ the spaces $\mathbb{H}_{W_\beta}({\D})$ satisfy the three conditions (i), (ii), and (iii).  By invoking Theorem \ref{Theorem 1} the $n$-best approximation problems have solutions in those spaces.

\begin{remark} (For the spaces $\mathbb{H}_{W_\beta}({\D}), \beta>1$)
The spaces $\mathbb{H}_{W_\beta}({\D}), \beta>1,$ do not fall into the category governed by Theorem \ref{Theorem 1}, as, owing to (\ref{owing to}), the condition (\ref{infinite norm}) is not satisfied. Recall the Sobolev Embedding theorem asserting that  $W^{k,p}\subset C^{r,\alpha}$ if $m<pk, \frac{1}{p}-\frac{k}{m}=-\frac{r+\alpha}{m}, m$ is the dimension. In our case $m=1, p=2, k=\frac{\beta}{2},$ and, in particular, $\beta =pk>1.$ It hence concludes that the functions in the spaces $\mathbb{H}_{W_\beta}({\D}), \beta>1,$ are all continuously extendable to the closed unit disc, and the norm square $\|k^\beta_a\|^2=k^\beta_a(a)$ given by (\ref{owing to}) does not have singularity for $|a|=1.$
  Based on these, as well as continuity of the inner product, we conclude existence of $n$-best approximations of the spaces for $\beta>1.$
\end{remark}

The results of this section are summarized as

\begin{theorem} For all Hardy-Sobolev spaces $\mathbb{H}_{W_\beta}({\D}), -\infty <\beta <\infty,$ there exist solutions to the $n$-best kernel approximation problem.
\end{theorem}

\section{Stochastic $n$-best Approximation}

Let $(\Omega, \mathcal{F}, d\mathbb{P})$ be a probability space, and, as in the previous sections, $\mathcal{H}$ be a RKHS of analytic functions in $\D$ with reproducing kernel $K_w, w\in \D.$ We will be studying stochastic signals $f(z,{\xi}),$ where $\xi\in \Omega$ and $z\in \D:$ We assume that for a.s. $\xi\in \Omega, f(\cdot, \xi)$ is a function in $\mathcal{H};$ and, for a.e. $z\in {\D}, f(z,\cdot)$ is a random variable. We will use the notation $f_\xi (z)=f(z,\xi).$ Associated with the probability space and the RKHS we define a Bochner type space (\cite{Boch})
\begin{eqnarray}
L^2(\mathcal{H},\Omega)=\{f: {\D}\times \Omega\to {\C}\ :\ \|f\|^2_{L^2(\mathcal{H},\Omega)}= E_\xi\|f_\xi\|^2_{\mathcal{H}}<\infty\},
\end{eqnarray}
where $E_\xi$ denotes the expectation, and precisely,
\[ E_\xi\|f_\xi\|^2_{\mathcal{H}}=\int_\Omega \|f_\xi\|^2_{\mathcal{H}}d\mathbb{P}(\xi).\]
We often use the simplified notation $\mathcal{N}=L^2(\mathcal{H},\Omega).$

The following theorem generalizes the existence result for stochastic $n$-best approximation for the Hardy space (see \cite{QQD}) to the RKHSs satisfying the conditions (i), (ii), and (iii) as assumed in Theorem \ref{Theorem 1}.

\begin{theorem}\label{Theorem 2} Let $\mathcal{H}$ be a RKHS of analytic functions in the unit disc satisfying the conditions (i), (ii), and (iii) in \S 1, and $\Omega$ a probability space. Let $\mathcal{N}= L^2(\mathcal{H},\Omega)$ be the associated Bochner type space as above defined. Let $f$ be any non-zero random signal. Then for any positive integer $n,$ either of the following two cases holds: (1) For some $1\leq m_1\leq n,$ there exists an $m_1$-tuple of constant parameters $(a_1,\cdots,a_{m_1})\in {\D}^{m_1}$ such that $f$ is identical with the orthogonal expansion
\begin{eqnarray}\label{expressible}
f(z,\xi)=\sum_{k=1}^{m_1}\langle f_\xi,E_k\rangle_{\mathcal{H}} E_k(z);
\end{eqnarray}
 or (2) There exists an $n$-tuple of constant parameters $(a_1,\cdots,a_{n})\in {\D}^{n}$ such that
\begin{eqnarray}\label{inf1} \|f-\sum_{k=1}^n \langle f_\xi,E_k\rangle_{\mathcal{H}} E_k\|_{\mathcal{N}}\end{eqnarray}
attains its positive infimum over all possible $n$-orthonormal systems $\{E_k\},$
where in both cases, $\{E_k\}_{k=1}^{m}$ is the orthonormal system generated by $(\tilde{K}_{a_1},\cdots,\tilde{K}_{a_{m}}), 1\leq m\leq n.$
\end{theorem}

The proof of Theorem 2.1 of \cite{QQD} for the stochastic Hardy space case cannot be directly adopted, for, in the present case no density argument based on the boundary value of the given function on $\partial \D$ is available. In \S 3 we established, by using a new technical method, the pointwise convergence result (\ref{rest}) which is crucial to in proving Theorem \ref{Theorem 2}. \\

\noindent {\bf Proof.}
In the proof of Theorem \ref{Theorem 1} we already show that, for each $\xi$ outside an event in $\Omega$ of probability zero, there holds uniformly in $a^{(l)}_1,a^{(l)}_2,...,a^{(l)}_{n-1}, l=1,2,\cdots,$ that
\[ \lim_{l\to \infty}|\langle f_\xi,E_{a^{(l)}_1a^{(l)}_2...a^{(l)}_n}\rangle_{\mathcal{H}}|^2=0.\]
By using the Cauchy-Schwarz inequality for the space $\mathcal{H}$ we have a dominating function of the function sequence on the left hand side:
\[|\langle f_\xi,E_{a^{(l)}_1a^{(l)}_2...a^{(l)}_n}\rangle_{\mathcal{H}}|^2
\leq \|f_\xi\|^2\in L^1(\Omega).\]
Then the Lebesgue dominated convergence theorem be invoked to conclude
\begin{eqnarray}\label{based on} \lim_{l\to\infty}E_\xi |\langle f_\xi,E_{a^{(l)}_1a^{(l)}_2...a^{(l)}_n}\rangle_{\mathcal{H}}|^2=0\end{eqnarray}
uniformly in $a^{(l)}_1,a^{(l)}_2,...,a^{(l)}_{n-1}, l=1,2,\cdots.$ Based on this the contradiction argument used in proving Theorem 2.1 of \cite{QQD} may be adopted to conclude the theorem.$\hfill\square$\\

\begin{remark}(Impact to Algorithm)
 This present paper only treats the existence aspect of the $n$-best problem. Based on the obtained estimates, however, a mathematical algorithm to actually get a solution is now on its way. We now cite the crucial step to reduce the problem to a global optimization one of a differential function defined in a compact set. Separate studies will be devoted to the computation aspect.
 To have an $n$-best approximation algorithm we are under the assumption that the given function $f$ is not expressible by any $m$-linear combination of multiple kernels for $m\le n-1.$ This implies that $d_f(n-1)>0,$ and there exists an $n$-tuple $(b_1,\cdots,b_n)\in {\bf D}^n$ such that for some $\epsilon>0$
\begin{eqnarray}\label{observation} \|f-P_{b_1\cdots b_n}f\|=d_f(n-1)-\epsilon.\end{eqnarray}
By using (\ref{rest}) one can find $\delta>0$ such that if
$|a_n|>1-\delta,$ then
\[\|f-P_{a_1\cdots a_{n-1}a_n}f\|>d_f(n-1)-\epsilon\]
for any $a_1,\cdots,a_{n-1}$ in $\D.$  Since $P_{a_1\cdots a_{n-1}a_n}f$ is symmetric in $a_1,\cdots,a_n,$ we conclude
\[\|f-P_{a_1\cdots a_{n-1}a_n}f\|>d_f(n-1)-\epsilon\]
whenever $|a_k|>1-\delta$ for some $k=1,\cdots,n.$
 Owing to the observation (\ref{observation}) we have $d_f(n)\leq d_f(n-1)-\epsilon.$ When the infimum attains at $(\tilde{a}_1,\cdots ,\tilde{a}_{n-1},\tilde{a}_n),$ that is,
 \[ \|f-P_{\tilde{a}_1\cdots \tilde{a}_{n-1}\tilde{a}_n}f\|=d_f(n),\]
 there holds $ |\tilde{a}_k|\leq 1-\delta$ for all $k=1,\cdots,n.$
 This concludes that the global minimum value $d_f(n)$ is only attainable in the compact set $ \overline{(1-\delta)\D}^n.$ In view of the relation (\ref{based on}) the same conclusion holds for the stochastic case.
\end{remark}

\section*{Acknowledgement}
The author wishes to express his sincere thanks to Y.B. Wang, K.H. Zhu and W. Qu for the previous works done with the author, as well as for their inspiriting discussions on this and related subjects.

\end{document}